\documentclass[11pt]{article}

\usepackage[latin1]{inputenc} % accents
\usepackage{latexsym}
\usepackage{graphicx}         % pour d'eventuelles figures
\usepackage{epsfig}           % (preferer graphicx, si possible)
\usepackage{amsmath}          % AMSTEX
\usepackage{amsfonts}
\usepackage{amssymb}
\usepackage{fancyhdr}

\DeclareGraphicsExtensions{.eps,.ps,.pdf,.bmp}

\newcommand{\titre}[1]{\begin{center}
                      {\Large{\bf #1}}
                      \end{center}}

\makeatletter
\renewcommand{\section}{\@startsection {section}{1}{0pt}
    {-3.5ex \@plus -1ex \@minus -.2ex}
    {2.3ex \@plus.2ex}
    {\normalfont\normalsize\bfseries}}
\newcommand{\SUMMARY}{\@startsection {section}{1}{0pt}
    {-3.5ex \@plus -1ex \@minus -.2ex}
    {2.3ex \@plus.2ex}
    {\normalfont\normalsize\centering\bfseries}*{{\bf R\'esum\'e}}}
    \makeatother

\begin{document}

%\pagestyle{fancy}
%\begin{center}
%\lhead{{\bf\hspace*{1.5cm}Troisième Congr\`es de la SM2A,
%Marrakech-Maroc, 10-13 Septembre 2012}}
%\end{center}
\rhead{}
%\lfoot{}
%\foot{(*) .}
%%%%%%%%%%%%%%%%%%%%%%%%%%%%%%%%%%%%%%%%%%%%%%%%%%%%%%%%%%%%%%%%%%%%%%%%%%%%
%               PRIERE DE NE RIEN MODIFIER AU-DESSUS                       %
%%%%%%%%%%%%%%%%%%%%%%%%%%%%%%%%%%%%%%%%%%%%%%%%%%%%%%%%%%%%%%%%%%%%%%%%%%%%
\newtheorem{theorem}{Theorem}[section]
\newtheorem{lemma}{Lemma}[section]
\newtheorem{definition}[theorem]{Definition}
\newtheorem{proposition}[theorem]{Proposition}
\newtheorem{corollary}[theorem]{Corollary}
\newtheorem{remark}{\bf Remark\/}
\newenvironment{prof}[1][Proof]{\noindent\textbf{#1.\ }}{\hfill~\rule{0.5em}{0.5em}}
%%%%%%%%%%%%%%%%%%%%%%%%%%%%%%%%%%%%%%%%%%%%%%%%%%%%%%%%%%%%%%%%%%%%%%%%%%%%
%                                                                          %
%        LES AUTEURS SONT PRIES DE FOURNIR  LES BONS ARGUMENTS             %
%       AUX COMMANDES titre, orateur, auteur, motsclefs CI-DESSOUS,        %
%                                                                          %
%%%%%%%%%%%%%%%%%%%%%%%%%%%%%%%%%%%%%%%%%%%%%%%%%%%%%%%%%%%%%%%%%%%%%%%%%%%%
%                                                                          %
% Remplacer "Titre" par le titre de votre comunication                     %
%                                                                          %
% Remplacer "Orateur" par le Nom de l'orateur                              %
%                                                                          %
% Remplacer "Auteurs" par les Noms des autres auteurs                      %
%                                                                          %
% Remplacer "Adresse 1 & 2" par l'adresse des auteurs                      %
%                                                                          %
% Remplacer "mail" par l'adresse électronique de correspondance            %
%                                                                          %
%%%%%%%%%%%%%%%%%%%%%%%%%%%%%%%%%%%%%%%%%%%%%%%%%%%%%%%%%%%%%%%%%%%%%%%%%%%%
\vspace*{-0.5cm} \titre{Bahadur efficiency of nonparametric test for
independence based on $L_1$-error}\underline{Noureddine Berrahou} $^{(1)*}$
and Lahcen Douge $^{(1)}$
\medskip

\indent$^1$ {\footnotesize \it  FSTG, Universit\'e Cadi Ayyad,
B.P. 549 Marrakech Maroc\\

$^*$Corresponding author : n.berrahou@uca.ma, lahcen.douge@uca.ma. }

\medskip
\vspace*{+0.2cm}
%%%%%%%%%%%%%%%%%%%%%%%%%%%%%%%%%%%%%%%%%%%%%%%%%%%%%%%%%%%%%%%%%%%%%%%%%%%%
%                                                                          %
%                   DEBUT DE LA COMMUNICATION                              %
%                                                                          %
%%%%%%%%%%%%%%%%%%%%%%%%%%%%%%%%%%%%%%%%%%%%%%%%%%%%%%%%%%%%%%%%%%%%%%%%%%%%
{\bf Abstract:} We introduce new test statistic to test the
independence  of two multi-dimensional random variables. Based on
the $L_1$-distance and the historgram density estimation method, the
test is compared via Bahadur relative efficiency to several tests
available in the literature. It arises that our test reaches better
performances than a number of usual tests among whom we cite the
Kolmogorov-Smirnov test. Beforehand, large deviation result is
stated for the associated statistic. The local asymptotic optimality
relative to the test is also studied.


\section{Introduction}
Consider a sample of $\mathbb{R}^d \times\mathbb{R}^{d^{'}}$-valued
random vectors $(X_1,Y_1),\ldots,(X_n,Y_n)$ with independent and
identically distributed (i.i.d.) pairs defined on the same
probability space. The density function of $(X,Y)$ is denoted by
$f$, while $f_1$ and $f_2$ stand for the density function of $X$ and
$Y$ with respect to the Lebesgue measure $\lambda_1$ and
$\lambda_2$, respectively. In this paper, we are concerned with the
problem of independence testing and the comparison of performances
of various tests. Here, we propose to test the  hypothesis that $X$
and $Y$ are independent on the basis of the shape of its density
from the $L_1$ point of view. More precisely, we consider the test
of the following null hypothesis
$$
H :\int\int|f(x,y) - f_1(x)f_2(y)|dxdy = 0.
$$
%against the alternative
%$$A :
%\int\int|f(x,y) - f_1(x)f_2(y)|dxdy >0.
% $$
For this purpose, we consider the histogram density estimates of the
density of $(X, Y)$, $X$ and $Y$. In order to define its, let $P_n
=\{A_{n,j}, j\geq1\}$ and $Q_n =\{B_{n,j}, j\geq1\}$  be a cubic
partitions of $\mathbb{R}^d$ and $\mathbb{R}^{d^{'}}$ respectively,
and denote by $\nu_n$, $\mu_{n,1}$ and $\mu_{n,2}$ the empirical
measures associated with the samples $(X_1,Y_1),\ldots,(X_n,Y_n)$,
$X_1,\ldots,X_n$, and $Y_1,\ldots,Y_n$, respectively. The histogram
density estimators of $f$ and $f_1$  are respectively defined by
\begin{eqnarray*}
f_{n}(x,y)=\frac{\mu_n(A_{n,j}\times
B_{n,j})}{\lambda_1(A_{n,j})\lambda_2(B_{n,j})} \quad \
\mbox{whenever}\ (x, y)\in A_{n,j}\times B_{n,j}
\end{eqnarray*}
and
\begin{eqnarray*}
f_{n,1}(x)=\frac{\mu_{n,1}(A_{n,j})}{\lambda_1(A_{n,j})} \,
\mbox{whenever}\ x\in A_{n,j},
\end{eqnarray*}
with $f_{n,2}$ is defined similarly as $f_{n,1}$. To perform the test of the
hypothesis $H$, we consider the following statistic
\begin{eqnarray*}
V_{n}=\int\int |f_{n}(x,y)-f_{n,1}(x)f_{n,2}(y)|dxdy.
\end{eqnarray*}

% The rejection regions relative to the test is given by
%$$R=\{V_n>C\},$$
%for some positive threshold $C$.
The work performed here deals with the $L_1$-distance large
deviation result for statistic $V_{n}$ and the comparison of their
performances in the Bahadur efficiency sense with some other tests
available in the literature. Furthermore, the local asymptotic
optimality relative to the test is studied. The independence
hypothesis may be expressed in various ways. Either the case
pertaining with the density and the $L_1$-distance in study in
this paper, the independence of $X$ and $Y$, based on the equality
of the joint distribution function and the product of its
marginals, may be defined by the following hypothesis
$$
H_{1}: F(x,y)=F_1(x)F_2(y),\quad \mbox{for all}\quad (x,y)\in
\mathbb{R}^d\times \mathbb{R}^{d^{'}},
$$
where $F$ is the distribution function of $(X,Y)$, and $F_1$, $F_2$
are the marginal distribution functions. A number of statistics have
been proposed to test the hypothesis $H_{1}$  or any other
hypothesis form against the general alternative that the equality in
$H_{1}$ is violated at least one point or narrower classes of
alternatives. Among the most studied ones in the univariate case, we
quote the Kolmogorov-type statistic has been introduced by
\cite{Blum}. It is defined by
$$
\Gamma_n=\sup_{-\infty<x,y<\infty}\mid F_n(x,y)-F_{n,1}(x)F_{n,2}(y)\mid,
$$
where $F_n$, $F_{n,1}$, and $F_{n,2}$  the empirical distribution
function, constructed from the initial sample and its components.
Properties of $\Gamma_n$ and its various versions have been
investigated in several papers among whom we quote
\cite{Deheuvels1982} and \cite{Nikitin1988}. On the basis of the
empirical distribution function several other tests have been
proposed for testing $H$. The well-known statistics are
\begin{eqnarray}\label{eqn1}
B^{k}_{n,q_1,q_2}&=&\int\int\left[ F_n(x,y)-F_{n,1}(x)F_{n,2}(y)\right]^{k}\nonumber\\
&&\hspace{3.5cm}\times q_1(F_{n,1}(x))q_2(F_{n,2}(y))
dF_{n,1}(x)dF_{n,2}(y),
\end{eqnarray}
where $k$ is a positive integer, and $q_1$ and $q_2$ are nonnegative
weight functions on $(0,1)$. The statistic $B^{1}_{n,q_1,q_2}$ for
$q_1(u)=q_2(u)=\sin(\pi u)$ has been introduced by Koziol and Nemec
(1979). For $k=2$ and $q_1=q_2=1$, a statistic equivalent to
(\ref{eqn1}) was proposed by \cite{Hoeffding}, with its properties
later being studied by \cite{Blum}, \cite{Cotterill} and
\cite{Nikitin1986}, and by \cite{DeWet} for weights $q_1$ and $q_2$
not equal to $1$.
%Performances of the independence tests based on the statistics $B^{1}_{n,q_1,q_2}$ and $B^{2}_{n,q_1,q_2}$  in the Bahadur sense have been compared to some usual tests in Nikitin (1995).
Similar to  the Durbin-type statistic in the goodness-of-fit
testing, \cite{Dmitrieva} proposed the statistic
$$
 M_n=\sup_{-\infty<x<\infty}\left| \int\left(F_n(x,y)-F_{n,1}(x)F_{n,2}(y)\right) dF_{n,2}(y)\right|.
$$
The author also investigated large deviations of this statistic and
their various versions under the independence hypothesis. As to
linear rank statistics for testing $H_1$, they usually are written
in the form
$$
 T_n=n^{-1}\sum_{i=1}^{n}a_{n,1}(R_{i}/(n+1))a_{n,2}(S_{i}/(n+1)),
$$
where $R_{i}$ is the rank of $X_i$ among $X_1,\ldots,X_n$,  $S_{i}$
is the rank of $Y_i$ among $Y_1,\ldots,Y_n$, and $a_{n,1}$,
$a_{n,2}$ are some real functions on $[ 0, 1]$. As an usual test
very close to the linear rank test, we cite the  Kendall rank
correlation coefficient $\tau_n$ (see \cite{Hettmansperge}). It is
defined by
$$
\tau_n=\frac{1}{n(n-1)}\sum_{i\neq j}\mbox{sign}(R_i-R_j)\mbox{sign}(S_i-S_j).
$$
Large deviations and exact slope of this statistic have been
investigated by \cite{sievers} and \cite{woodworth}. We refer also
to \cite{Nikitin1995} where the comparison of these statistics by
their asymptotic efficiency is done. A statistic very close to $V_n$
was proposed by \cite{gretton}. It is defined by
\begin{eqnarray*}
  L_n(\nu_n,\mu_{n,1}\times \mu_{n,2})=\sum_{A\in P^{'}_{n}}\sum_{B\in Q^{'}_{n}}\mid\nu_n(A\times B)-\mu_{n,1}(A)\mu_{n,2}( B)\mid,
\end{eqnarray*}
where $P^{'}_{n}$ and $Q^{'}_{n}$ are finite partitions of
$\mathbb{R}^d$ and $\mathbb{R}^{d^{'}}$ respectively. Asymptotic
properties of the statistic $L_n$ and, particular, corresponding
large deviations were investigated in  \cite{gretton}. However, the
methods applied in this later paper are not sufficiently strong  to derive
the Bahadur efficiency. Under the independence hypothesis, these
authors proved only that for all $0<\epsilon_1, 0<\epsilon_2$ and
$0<\epsilon_3$,
$$
\mathbf{P}\left(L_n(\nu_n,\mu_{n,1}\times \mu_{n,2})>\epsilon_1
+\epsilon_2+\epsilon_3\right)\leq 2^{m_n m^{'}_n}e^{-n\epsilon_1^{2}/2}
+2^{m_n}e^{-n\epsilon_2^{2}/2}+2^{m^{'}_n}e^{-n\epsilon_3^{2}/2}.
$$
A significantly more exact large deviation result  is given by the
following Theorem, proved in this work.
\begin{theorem}
Under the independence hypothesis $H$ one has
$$
 \lim_{n\rightarrow\infty}n^{-1}\log \mathbf{P}\left(V_n>\lambda\right)=-g(\lambda)=
 -\frac{\lambda^{2}}{2}(1+o(1)), \quad \lambda\rightarrow0,
$$
where the function $g$ is continuous for small $\lambda>0$.
\end{theorem}

%\medskip
%{\bf Conclusion:}

%%%%%%%%%%%%%%%%%%%%%%%%%%%%%%%%%%%%%%%%%%%%%%%%%%%%%%%%%%%%%%%%%%%%%%%%%%%%
%               REFERENCES BIBLIOGRAPHIQUES                                %
%%%%%%%%%%%%%%%%%%%%%%%%%%%%%%%%%%%%%%%%%%%%%%%%%%%%%%%%%%%%%%%%%%%%%%%%%%%%
\bibliographystyle{plain}

\begin{thebibliography}{99}

\bibitem{Blum}
{\sc Blum, J.~R., Kiefer, J., and Rosenblatt, M. (1961)}.
\newblock {\em Distribution free tests of independence based on the sample
  distribution function},
\newblock Ann. Math. Statist. 32, 485--498.


\bibitem{Cotterill}
{\sc Cotterill, D.~S. and Cs{\"o}rg{\H{o}}, M. (1985)},
\newblock {\em On the limiting distribution of and critical values for the
  {H}oeffding, {B}lum, {K}iefer, {R}osenblatt independence criterion},
\newblock Statist. Decisions, 3 (1-2), 1--48.

\bibitem{DeWet}
{\sc De~Wet, T. (1980)},
\newblock {\em Cram\'er-von {M}ises tests for independence},
\newblock J. Multivariate Anal., 10 (1), 38--50.

\bibitem{Deheuvels1982}
{\sc Deheuvels, P. (1982)},
\newblock {\em Some applications of the dependence functions to statistical
  inference: nonparametric estimates of extreme values distributions, and a
  {K}iefer type universal bound for the uniform test of
  independence},
\newblock In: Nonparametric statistical inference, {V}ol. {I}, {II}
  ({B}udapest, 1980), volume~32 of {\em Colloq. Math. Soc. J\'anos Bolyai},
  pages 183--201. North-Holland, Amsterdam.


\bibitem{Dmitrieva}
{\sc Dmitrieva, A.~G. (1988)},
\newblock {\em Asymptotic behavior of a statistic of {D}urbin type for a test for
  independence},
\newblock Vestnik Leningrad. Univ. Mat. Mekh. Astronom., (vyp. 1),
  102--104, 119.

\bibitem{gretton}
{\sc Gretton, A. and Gy{\"o}rfi, L. (2010)},
\newblock {\em Consistent nonparametric tests of independence},
\newblock J. Mach. Learn. Res., 11, 1391--1423.

\bibitem{Hettmansperge}
{\sc Hettmansperger, T.~P. (1984)},
\newblock {\em Statistical inference based on ranks},
\newblock Wiley Series in Probability and Mathematical Statistics: Probability
  and Mathematical Statistics. John Wiley \& Sons Inc., New York.

\bibitem{Hoeffding}
{\sc Hoeffding, W. (1948)},
\newblock {\em A non-parametric test of independence},
\newblock Ann. Math. Statistics, 19, 546--557.

\bibitem{Nikitin1995}
{\sc Nikitin, Y. (1995)}.
\newblock {\em Asymptotic efficiency of nonparametric tests}.
\newblock Cambridge University Press, Cambridge.

\bibitem{Nikitin1986}
{\sc Nikitin, Y.~Y. (1986)},
\newblock {\em Large deviations and asymptotic efficiency of integral statistics for
  a test of independence},
\newblock In: Probability distributions and mathematical statistics
  ({R}ussian) ({F}ergana, 1983), pages 388--406, 497--498. ``Fan'', Tashkent.

\bibitem{Nikitin1988}
{\sc Nikitin, Y.~Y. and Pankrashova, A.~G. (1988)},
\newblock {\em Bahadur efficiency and local asymptotic optimality of some
  nonparametric independence tests},
\newblock Zap. Nauchn. Sem. Leningrad. Otdel. Mat. Inst. Steklov. (LOMI)
  166 (Issled. po Mat. Statist. 8), 8, 112--128, 189.

\bibitem{sievers}
{\sc Sievers, G.~L. (1969)},
\newblock {\em On the probability of large deviations and exact slopes},
\newblock Ann. Math. Statist. 40, 1908--1921.

\bibitem{woodworth}
{\sc Woodworth, G.~G. (1970)},
\newblock {\em Large deviations and {B}ahadur efficiency of linear rank statistics},
\newblock Ann. Math. Statist.  41, 251--283.


\end{thebibliography}

%%%%%%%%%%%%%%%%%%%%%%%%%%%%%%%%%%%%%%%%%%%%%%%%%%%%%%%%%%%%%%%%%%%%%%%%%%%

\end{document}